\documentclass[letterpaper,10pt,conference]{article}
\usepackage{url}
\usepackage{amssymb}
\usepackage{graphicx}
\usepackage{amsfonts}
\usepackage{amsbsy}
\usepackage{fancybox}
\usepackage{float}
\usepackage{latexsym}
\usepackage{makeidx}
\usepackage{psfrag}
\usepackage{lscape}

\usepackage{times}
\usepackage[intlimits]{amsmath}
\usepackage{hhline}

\floatstyle{ruled}
\newfloat{matlab}{!!!!htbp}{mat}[section]
\floatname{matlab}{Script Matlab}

\newtheorem{definition}{Definition}

\newtheorem{assumption}{Assumptions}
\newtheorem{remark}{Remark}

\newcommand{\myfig}[2]{\centerline{\includegraphics[width=#1]{#2}}}
\newcommand{\bequ}{\begin{equation} }
\newcommand{\eequ}{\end{equation}}

\newcommand{\bary}{\begin{array}}
\newcommand{\eary}{\end{array}}
\newcommand{\blem}{\begin{lem} }
\newcommand{\elem}{ \end{lem} }
\newcommand{\bpro}{\begin{prop} }
\newcommand{\epro}{ \end{prop} }
\newcommand{\bass}{\begin{assumption} }
\newcommand{\eass}{ \end{assumption} }
\newcommand{\bthm}{\begin{thm} }
\newcommand{\ethm}{ \end{thm} }
\newcommand{\bpb}{\begin{prob} }
\newcommand{\epb}{ \end{prob} }
\newcommand{\bdif}{\begin{definition} }
\newcommand{\edif}{ \end{definition}}
\newcommand{\bcor}{\begin{cor} }
\newcommand{\ecor}{ \end{cor} }
\newcommand{\brem}{\begin{remark} }
\newcommand{\erem}{ \end{remark} }
\newcommand{\bmatlab}{\begin{matlab}}
\newcommand{\ematlab}{\end{matlab}}

\newcommand{\bmat}[1]{\left [ \bary{#1}}
\newcommand{\emat}{\eary \right ] }
\newcommand{\bvec}[1]{\left ( \bary{#1}}
\newcommand{\evec}{\eary \right ) }

\newcommand{\1}{{{I}}}
\newcommand{\0}{{\boldsymbol{0}}}

\newcommand{\be}{\begin{equation}}
\newcommand{\ee}{\end{equation}}

\def\CheckPDFoutput{%
\CheckPDFoutput%
\ifx\unprotect\undefined%
 \DeclareGraphicsRule{.jpg}{bmp}{}{}%
  \else%
  \pdfoutput=1%
\fi%
}

\usepackage{epstopdf}
\usepackage{fancyhdr}
\usepackage{color}
\usepackage{booktabs}

\title{ Randomized and robust methods for uncertain systems using \RR, with applications to DEMETER satellite benchmark}

\author{Mohammadreza Chamanbaz\footnote{Corresponding author: mrchamanbaz@gmail.com.},\\ Singapore University of Technology and Design, Singapore \\[1em]
Fabrizio Dabbene, CNR-IEIIT Politecnico di Torino, Italy\\[1em]
Dimitri Peaucelle, \\LAAS-CNRS, Universit\'e de Toulouse, CNRS, Toulouse, France\\[1em]
Christelle Pittet, CNES, Toulouse, France\\[1em]
Roberto Tempo, CNR-IEIIT Politecnico di Torino, Italy}

\date{{\it {\bf Aerospace Lab} \\ Special issue on aerial robotics. \\ \vspace*{2mm} {\footnotesize Final version, \today.}}}

\newcommand{\RR}{\textsc{R-RoMulOC}}

\begin{document}
\maketitle
\thispagestyle{empty}
\pagestyle{empty}

\begin{abstract}
\RR~ is a freely distributed toolbox which aims at making easily available to the users different optimization-based methods for  dealing with uncertain systems. It implements both deterministic LMI-based results, that provide guaranteed performances for all values of the uncertainties, and probabilistic randomization-based approaches, that guarantee performances for all values of the uncertainties except for a subset with arbitrary small probability measure. The paper is devoted to the description of these two approaches for analysis and control design when applied to a satellite benchmark proposed by CNES, the French Space Agency. The paper also describes the modeling of the DEMETER satellite and its integration into the \RR~ toolbox as a challenging test example. Design of state-feedback controllers and closed-loop performance analysis are carried out with the randomized and robust methods available in the \RR~ toolbox. 
\end{abstract}
\section{Introduction}

The last decades have witnessed an increase of interest in the area of analysis and design of systems in the presence of uncertainty. This is due to the continuous development of novel and efficient theoretical and numerical tools for robustness  (ability of the system to maintain stability and performance under large variations of the system parameters), see \cite{petersen_robust_2014} for a recent overview.

In particular, two main paradigmatic approaches have gained popularity. On one side, the worst-case, or deterministic, paradigm aims at guaranteeing a desired  level of performance \textit{for all} system's configurations. This approach has largely benefited from the introduction of the linear matrix inequalities (LMIs)  formalism, which led to many important results, allowing to tackle a large variety of uncertainty models and performances requirements. Recently, the corresponding numerical tools have been collected in a \textsc{Matlab} toolbox named Robust Multi Objective Control toolbox (\textsc{RoMulOC}) \cite{romulocconf}.
The toolbox provides different functions for describing and manipulating uncertain systems, and for building LMI optimization problems related to robust multiobjective control problems. We refer to \cite{petersen_robust_2014} for and extensive review of deterministic and probabilistic methods in robust control design and analysis.  

The deterministic  approach can be seen as ``pessimistic," in the sense that the guaranteed (and certified) performance is usually significantly worse than the actual worst case performance, due to unavoidable conservatism of the developed methodologies. This fact motivated the
introduction of a probabilistic approach \cite{tempo_randomized_2012,calafiore_research_2011}, which consists in testing a finite number of configurations among the infinitely many admissible ones. This approach is said to be ``optimistic," in the sense that even if a level of performance is valid for all tested cases, it may not hold for some of the unseen instances. However, rigorous theoretical results, based on large-deviation inequalities, have been derived to bound the probability of performance violation. This theory has now reached a good level of maturity, and the main algorithms have been coded in the Randomized Algorithms Control Toolbox (RACT) \cite{tremba_ract:_2008} which can be freely downloaded from \url{http://ract.sourceforge.net/pmwiki/pmwiki.php/}. This toolbox allows the user to define and manipulate various types of probabilistic uncertainties, providing efficient sampling algorithms for the different uncertainty types commonly encountered in robust control. Furthermore, it includes sequential and batch randomized algorithms for control systems design.

It is important to remark that these two paradigms are not in competition, but they represent complementary approaches that provide additional tools to the systems engineer for the design of control system under uncertainty. Inspired by these considerations,  a joint effort between the two teams at the core of \textsc{RoMulOC} and RACT has been recently carried out, with the aim of merging the features of the two toolboxes in an integrated framework. This lead to the development of \RR. The main feature of this toolbox is to allow the user to input the system's description only once, using the well tested formalism of \textsc{RoMulOC}. Then, both deterministic and probabilistic methods can be applied on the same system, efficiently moving from a deterministic to a probabilistic description of the uncertainty, by simply changing some parameters in the code.

As the two tools from which it originates, \RR~is freely distributed, and can be downloaded at \url{http://projects.laas.fr/OLOCEP/rromuloc/}. We refer the interested reader to this webpage for a detailed list of references to the various worst-case and probabilistic methods which are coded in \RR. For a description of the \RR~toolbox, the reader is referred to \cite{rromulocconf}.

In this paper,  the effectiveness of the toolbox is shown by introducing the modeling of the DEMETER satellite \cite{pit:arz/rocond06} in \RR~toolbox. Then, we show how the design of state-feedback controllers and analysis of closed-loop performance can be performed with the randomized and robust methods available in the \RR~toolbox.\\

\noindent
{\bf Notation\\}
$\1_n$ stands for the identity matrix of dimension $n$.
$A^T$ is the transpose of $A$. $\{A\}^{\mathcal S}$ represents the symmetric matrix $\{A\}^{\mathcal S}=A+A^T$.
${\sf Tr}(A)$ is the trace of $A$.
$A\succ (\succeq)B$ means $A-B$ is positive (semi-)definite.
${\sf diag}\bmat{ccc} \cdots  & F_i & \cdots \emat$ is a block-diagonal matrix whose diagonal blocks are $F_i$.
The symbol $\otimes$ refers to Kronecker product. Given vectors $v,w\in \mathbb R^3$, the matrix $v^\times\in\mathbb R^{3\times 3}$ is a skew-symmetric matrix defined such that $v\times w= v^\times w$, i.e. 
\[v^\times=\bmat{ccc}0&-v_z&v_y\\v_z&0&-v_x\\-v_y&v_x&0\emat\]
for $v= \bmat{ccc}v_x & v_y & v_z\emat^T$. The three-dimensional sphere $\mathbb S^3$ is parameterized by quartenions $q\in\mathbb R^4$ satisfying the constraint $|q| = 1$. 
Finally, star-product describes Linear-Fractional Transformations (LFT)
\[M_a+M_b\Delta(\1-M_d\Delta)^{-1}M_c=\Delta\star\bmat{c|c}M_d&M_c\\\hline M_b&M_a\emat.\]

\section{DEMETER benchmark}

\begin{figure}
\centering
\includegraphics[width=0.5\columnwidth]{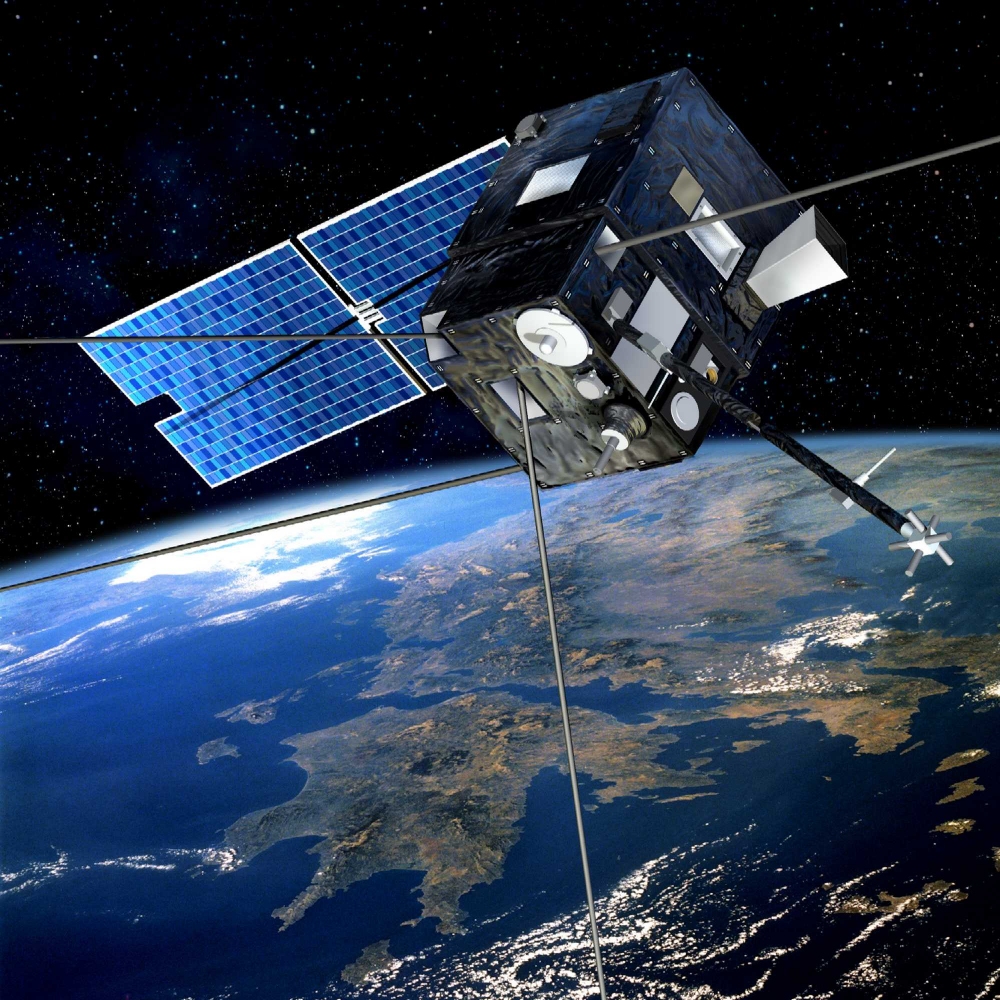}
\caption{ DEMETER satellite. \textcopyright CNES November 2003, ill. D. Ducros.  }
\label{fig: demeter}
\end{figure}

DEMETER is a satellite of the CNES Myriade series. Launched in 2004, it observed electric and magnetic signals in Earth's ionosphere for more than 6 years. Its characteristic is to be composed of a central body and four long and flexible appendices---as shown in Figure \ref{fig: demeter}---oriented in different directions and fixed to the rigid-body at different positions distinct from the center of gravity. The model of this satellite has been provided as a benchmark in \cite{pit:arz/rocond06}. This model with uncertainties is revisited in the following. A specific function incorporated in \RR~allows to generate variants of the complete benchmark. The variants are such that the user can generate models of various sizes, both in terms of order of the plant and in terms of the number of uncertainties involved. 
 
\subsection{Nonlinear model without flexible modes}

Assuming full actuation for attitude control $u\in\mathbb{R}^3$ and modeling in the body-fixed frame, the nonlinear dynamics of the satellite are 
\begin{equation}\label{eq:nonlinear model}
J\dot\omega+\omega^\times J\omega=u
~~~~,~~~~
\dot q=\frac{1}{2}\bmat{cc} -\omega^\times & \omega\\-\omega^T &0\emat q,
\end{equation}
where $\omega\in\mathbb{R}^3$ is the rotational velocity of the satellite body-fixed frame with respect to the inertial frame, $J\in\mathbb{R}^{3\times 3}$ is the symmetric positive definite matrix corresponding to its moment of inertia and $q\in\mathbb{S}^3$ are the quaternion coordinates.
A classical control problem related to this nonlinear model is to build an ideal state-feedback control law $u^\star(\omega,q)$ guaranteeing  global stability. A more involved problem is to take into account in the design phase implementation issues such as saturation of reaction wheels, sensor delays and failures, periodic sub-actuated character of magneto-torquers, etc. The model complexity depends on the considered actuators. For example, considering reaction wheel control, the model becomes
\begin{equation}\label{eq:nonlinear model with reaction wheels}
J\dot\omega+\omega^\times (J\omega+h)=-T+T_{{\sf ext}}
~~~,~~~
\dot h=T
~~~,~~~
\dot q=\frac{1}{2}\bmat{cc} -\omega^\times & \omega\\-\omega^T &0\emat q,
\end{equation}
where $h\in\mathbb R^3$ is the vector of the angular {momenta} of the wheels and $T$ is the vector of the torques applied to the wheels and $T_ {{\sf ext}}$ is the external disturbances controller should reject.

\subsection{Linear model with flexible modes}

Let $\theta\in\mathbb{R}^3$ be the three axes angular deviation of the satellite from some reference constant orientation. The linearized  model of \eqref{eq:nonlinear model} is
\begin{equation}\label{eq: linearized model}
J\ddot\theta=u,
\end{equation}
which is a three dimensional double integrator. We remark that so far we assumed that the satellite is composed only of a rigid body. Unfortunately, this is not the case because of solar panels and other scientific equipment on board. At small pointing errors (the attitude control is required to have less than 0.1 degree precision) the flexibility of appendices is not negligible and needs to be considered in the model. The linearized model including flexible modes is \cite{pit:arz/rocond06} 
\begin{equation}\label{eq: linearized model with flexible modes}
\bmat{cc} J & J^{1/2}L \\ L^TJ^{1/2} & \1 \emat
\bvec{c} \ddot\theta \\\ddot\eta\evec
+
\bmat{cc} \0&\0\\2Z\Omega&\Omega^2\emat
\bvec{c} \dot\eta\\\eta\evec
=
\bmat{c} \1 \\\0\emat
u,
\end{equation}
where $\eta\in\mathbb R^{2n_f}$ is the vector of angular deviations in torsion and bending of the flexible appendices (up to $n_f=4$ in the DEMETER model), $L$ is a matrix modeling the cross influence of flexible modes on the rigid body which depends on how the appendices are attached to the rigid body,
$Z= {\sf diag}\bmat{ccc}\cdots &\zeta_iI_2&\cdots\emat$ is a diagonal matrix of all flexible mode damping factors and $\Omega={\sf diag}\bmat{ccc}\cdots &\omega_iI_2&\cdots\emat$ is a diagonal matrix of all flexible mode natural frequencies (the low damped oscillatory flexible dynamics are such that $\ddot\eta_i+2\zeta_i\omega_i\dot\eta_i+\omega_i^2\eta_i=-L_i^TJ^{1/2}\ddot\theta$). The same parameters apply for the bending and torsion effects, and in most cases one can assume that the appendices are identical ($\zeta_i=\zeta\,\forall i=1,\dots ,n_f$ and $\omega_i=\omega\,\forall i=,1\dots, n_f$). In \eqref{eq: linearized model with flexible modes}, the force $L^TJ^{1/2}\ddot\theta$ that acts on the flexible modes comes from the derivative of the angular momentum of the rigid body, and its symmetric feedback reaction on the rigid body is $J^{1/2}L\ddot\eta$. An analysis in frequency domain shows that only the first flexible modes of the appendices have significant influence on the system dynamics while all other flexible modes, including those of the solar panels, can be neglected. 

\subsection{Parametric uncertainties}

In \eqref{eq: linearized model with flexible modes} the matrix $L$, which is only due to positioning of the appendices, is assumed to be perfectly known. 
All other parameters, {\it i.e. }$J$, $\zeta_i$ and $\omega_i$ cannot be precisely measured on the earth due to gravity, hence are considered to be uncertain.
Damping ratio and natural frequencies $\zeta_i,\omega_i$  describe the first flexible modes of the four appendices. These appendices are of same length and same material and hence their flexible modes are almost identical. Yet there are discrepancies from one appendix to another which are not known. Damping ratio and natural frequencies are assumed to be bounded in the intervals
\[
\omega_i\in[~0.2\cdot 2\pi~,~0.6\cdot 2\pi~]
~~~,~~~
{\zeta_i}\in[~5\cdot 10^{-4}~,~5\cdot10^{-3}~]
~~~\forall i=1,\dots, 4.
\]
The inertia $J$ has the following nominal value  on ground 
\[
J_o=
\bmat{ccc}
J_{o11} & J_{o12} & J_{o13}\\
J_{o12}&J_{o22} & J_{o23}\\
J_{o13}&J_{o23}& J_{o33}
\emat
=
\bmat{ccc}
31.38 & -1.11 & -0.26\\
-1.11&21.19 & -0.78\\
-0.26&-0.78& 35.70
\emat.
\]
Uncertainties in $J$ are assumed to be at most of 30\% on the diagonal entries and $\pm 3$ on the off-diagonal entries. That is, for example, $J_{11}\in[~0.7 J_{o11}~,~1.3 J_{o11}~]=[~21.97~,~40.80~]$ and $J_{12}\in[~J_{o12}-3~,~J_{o12}+3~]=[~-4.11~,~1.89~]$.

\subsection{LFT modeling of uncertain matrices}

We first derive the LFT model of the $\bmat{cc} 2Z\Omega&\Omega^2\emat$ matrix.
Note that the uncertain matrices $\Omega$ and $Z$ are defined as a nominal matrix with normalized discrepancies around the nominal value. Hence, one can write $\Omega$ as 
\[
\Omega
=\omega_a\1+\omega_b\delta_\Omega
=\delta_\Omega\star\bmat{c|c}\0&\1\\\hline\omega_b\1&\omega_a\1\emat~~~\delta_\Omega={\sf diag}\bmat{cccc} \delta_{\omega_1}\1_2 & \delta_{\omega_2}\1_2 & \delta_{\omega_3}\1_2 & \delta_{\omega_4}\1_2 \emat,
\]
where $\omega_a=\frac{1}{2}(0.6\cdot 2\pi+0.2\cdot 2\pi)=0.4\cdot 2\pi$ is the mean between the two extreme values,
$\omega_b=\frac{1}{2}(0.6\cdot 2\pi-0.2\cdot 2\pi)=0.2\cdot 2\pi$ is the maximal deviation and $|\delta_{\omega_i}|\leq 1,~i = 1,\ldots,4$ are norm bounded uncertainties. The uncertain matrix $Z$ can be derive in a similar way
\[
Z
=\zeta_a\1+\zeta_b\delta_Z
=\delta_Z\star\bmat{c|c}\0&\1\\\hline\zeta_b\1&\zeta_a\1\emat,~~~\delta_Z={\sf diag}\bmat{cccc} \delta_{\zeta_1}\1_2 & \delta_{\zeta_2}\1_2 & \delta_{\zeta_3}\1_2 & \delta_{\zeta_4}\1_2 \emat,
\]
with $\zeta_a=\frac{1}{2}(5\cdot10^{-3}+5\cdot 10^{-4})=2.75\cdot 10^{-3}$ being the mean between the two extreme values,
$\zeta_b=\frac{1}{2}(5\cdot10^{-3}-5\cdot 10^{-4})=2.25\cdot 10^{-3}$  being the maximal deviation and $|\delta_{\zeta_i}|\leq 1,~i=1,\ldots,4$ are the norm bounded uncertainties. 
Using properties of the star-product we have
\[
\bmat{cc}2Z&\Omega\emat
=
\bmat{cc} \delta_{Z}&\0\\\0&\delta_{\Omega}\emat
\star
\bmat{cc|cc} \0&\0 & \1 &\0\\\0 &\0 &\0 &\1\\\hline 2\zeta_b\1&\omega_b\1 & 2\zeta_a\1&\omega_a\1\emat,
\]
and
\[\bary{l}
\bmat{cc}2Z\Omega&\Omega^2\emat
=\Omega\bmat{cc}2Z&\Omega\emat
\\
=
\bmat{ccc} \delta_{\Omega}&\0&\0\\\0&\delta_{Z}&0\\\0&\0&\delta_{\Omega}\emat
\star
\bmat{ccc|cc}
	\0  & 2\zeta_b\1&\omega_b\1 &  2\zeta_a\1&\omega_a\1\\
	\0 & \0 & \0 &\1&\0\\
	\0 & \0 & \0& \0 &\1\\\hline
	\omega_b\1 & 2\omega_a\zeta_b\1 &\omega_a\omega_b\1 &  2\omega_a\zeta_a\1&\omega_a^2\1
\emat
\eary.\]
We remark that the LFT defined in this way is minimal. An alternative is to build separately the LFTs for $2Z\Omega$ and $\Omega^2$ matrices and then to concatenate the two. This alternative gives an LFT with $\delta_\Omega$ repeated 3 times, which is clearly non-minimal.

\noindent
We next focus on the LFT modeling of the matrix depending on the uncertain matrix $J$. The difficulty can be observed arising from modeling  square-root of $J$. In \cite{pit:arz/rocond06} it is implicitly assumed that off-diagonal terms in $J$ are sufficiently small to be neglected in the computation of $J^{1/2}$. That is,  defining
\[
J=J_1+J_1^T+J_2
~~:~~
J_1=\bmat{ccc}0&J_{12}&J_{13}\\0&0&J_{23}\\0&0&0\emat
,~
J_2={\sf diag}\bmat{ccc}J_{11} & J_{22} & J_{33} \emat,
\]
it is assumed that $J^{1/2}\simeq J_2^{1/2}$. Then, to further simplify the model, the paper \cite{pit:arz/rocond06} makes the second assumption that the square root can be replaced by a first order approximation  $(J_{2a}+J_{2b}\delta_{J_2})^{1/2}\simeq J_{2a}^{1/2}+\frac{1}{2}J_{2b}\delta_{J_2}$. The relative error of this last approximation is less than 2\%, which is indeed reasonable. Based on this approximation, the minimal LFT model is such that $\delta_{J_2}$ is repeated twice. As we will show next, there is no reason for performing the first order approximation, and this can be avoided without increasing the size of the LFT.

\noindent
Two ways for improving the square root LFT modeling are explored next. The first one still assumes that $J^{1/2}\simeq J_2^{1/2}$ but avoids the first order approximation of the square root. To this end, define the following LFT modeling of the square root of inertias diagonal components
\[
J_2^{1/2}
=\hat J_{2a}+\hat J_{2b}\delta_{\hat J_2}
=\delta_{\hat J_2}\star\bmat{c|c}\0&\1\\\hline \hat J_{2b}&\hat J_{2a}\emat,
\]
where $\hat J_{2a}=\frac{1}{2}((1.3J_{2a})^{1/2}+(0.7J_{2a})^{1/2})$ is the mean between the two extreme values,
$\hat J_{2b}=\frac{1}{2}((1.3J_{2a})^{1/2}-(0.7J_{2a})^{1/2})$ is the maximal deviation, $\delta_{\hat J_2}={\sf diag}\bmat{ccc} \delta_{\hat J_{11}}&\delta_{\hat J_{22}} &\delta_{\hat J_{33}} \emat$ and $|\delta_{\hat J_{ii}}|\leq 1$ are the norm bounded uncertainties.
Using properties of the star-product one gets
\[
\bary{l}
\bmat{cc}J_2 & J_2^{1/2}L\\L^TJ_2^{1/2}&L^TL\emat
=
\bmat{c} J_2^{1/2} \\L^T\emat \bmat{cc} J_2^{1/2} & L\emat
\\=
\bmat{cc} \delta_{\hat J_2}&\0\\\0&\delta_{\hat J_2}\emat
\star
\bmat{cc|cc}
\0&\hat J_{2b}^2&\hat J_{2b}\hat J_{2a}&\hat J_{2b}L\\
\0&\0&\1&\0\\\hline
\1&\hat J_{2a}\hat J_{2b}&\hat J_{2a}^2&\hat J_{2a}L\\
\0&L^T\hat J_{2b}&L^T\hat J_{2a}&L^TL
\emat
\eary.
\]
Notice that---as in \cite{pit:arz/rocond06}---the uncertainties $\delta_{\hat J_2}$ are repeated only twice hence, the LFT size is not increased by precise modeling of the square root.

\noindent
Next, consider the cross inertia dependent matrix 
\[
J_1
=\hat J_{1a}+\hat J_{1b}\delta_{J_1}\hat J_{1c}
=\delta_{J_1}\star\bmat{c|c}\0&J_{1c}\\\hline J_{1b}& J_{1a}\emat,
\]
\[J_{1a}=\bmat{ccc} 0&J_{o12}&J_{o13}\\0&0&J_{o23}\\0&0&0\emat
,~
J_{1b}=\bmat{ccc} 3&3&0\\0&0&3\\0&0&0\emat
,~
J_{1c}=\bmat{ccc} 0&1&0\\0&0&1\\0&0&1\emat,
\]
\[
\delta_{J_1}={\sf diag}\bmat{ccc} \delta_{J_{12}}&\delta_{J_{13}} &\delta_{J_{23}} \emat
~~~:~~~
|\delta_{J_{ij}}|\leq 1.
\]
Using properties of the star-product we finally arrive at
\[\bary{l}
\bmat{cc} J_1+J_1^T+J_2&J_2^{1/2}L\\LJ_2^{1/2}&\1\emat
\\
=
{\sf diag}\bmat{c} \delta_{J_1} \\ \delta_{J_1} \\ \delta_{\hat J_2} \\ \delta_{\hat J_2} \emat
\star
\bmat{cccc|cc}
\0&\0&\0&\0&J_{1c}&\0\\
\0&\0&\0&\0&J_{1b}^T&\0\\
\0&\0&\0&\hat J_{2b}^2&\hat J_{2b}\hat J_{2a}&\hat J_{2b}L\\
\0&\0&\0&\0&\1&\0\\\hline
J_{1b}&J_{1c}^T&\1&\hat J_{2a}\hat J_{2b}&J_{1a}+J_{1a}^T+\hat J_{2a}^2&\hat J_{2a}L\\
\0&\0&\0&L^T\hat J_{2b}&L^T\hat J_{2a}&\1
\emat
\eary.\]
The second approach for improving the LFT modeling of the square-root $J^{1/2}$ needs first to question the relevance of modeling the coefficients of $J$ in intervals. The matrix $J$ is  symmetric positive definite which can be defined as $J=(J_o^{1/2}+\Delta_{\hat J})^2$ with an uncertain symmetric matrix $\Delta_{\hat J}$ constrained by a convex quadratic constraint
\[
X+Y\Delta_{\hat J}+\Delta_{\hat J}Y+\Delta_{\hat J}Z\Delta_{\hat J}\preceq\0 ~~,~~Z\succeq \1,
\]
where all $X$, $Y$ and $Z$ matrices are chosen symmetric to fit with the symmetric nature of $\Delta_{\hat J}$. The set also reads as
\[
(\Delta_{\hat J}-\Delta_o)Z(\Delta_{\hat J}-\Delta_o)\preceq \Delta_oZ\Delta_o -X,
\]
where $\Delta_o=-YZ^{-1}$ is the center of the set. Recall that $J=(J_o^{1/2}+\Delta_{\hat J})^2$ is (as formulated in \cite{pit:arz/rocond06}) a matrix whose $6$ independent coefficients are in intervals. The matrix $J$ can therefore be defined as the convex linear combination of $2^6$ vertices---denoted as $J^{[v]},~v=1,\ldots, 2^6$---and constructed taking all the extreme  combinations of the interval uncertainties. A natural way of defining  $X$, $Y$, $Z$ matrices is to impose the set to contain the convex combination of the square-roots of extremal values, that is the matrices $\Delta_{\hat J}^{[v]}={J^{[v]}}^{1/2}-J_o^{1/2}$
\bequ\label{e-elionvertices}
(\Delta_{\hat J}^{[v]}-\Delta_o)Z(\Delta_{\hat J}^{[v]}-\Delta_o)\preceq \Delta_oZ\Delta_o -X
~~
\forall v=1,\ldots, 2^6.
\eequ
A natural choice for the center of the set is to take the mean value of all vertices
\bequ\label{e-deltacenter}
\Delta_o=\frac{1}{2^6}\sum_{v=1}^{2^6} \Delta_{\hat J}^{[v]}.
\eequ
Of course, one aims at defining the smallest set containing the matrices $\Delta_{\hat J}^{[v]}$. It is rather easy to see that the size of the set is highly dependent on the matrix $\Delta_oZ\Delta_o -X$. The smaller it is, the smaller the set of $\Delta_{\hat J}$ matrices will be. It is suggested to minimize this matrix with respect to its Frobenius norm, which amounts to take 
\[
(X^*,Z^*)=\arg\min_{Z\succeq\1, (\ref{e-elionvertices})} {\sf Tr}(\Delta_oZ\Delta_o -X),
\]
and $Y^*=-\Delta_o{Z^*}^{-1}$. Having performed this LMI optimization, the inertia of the satellite is now defined as
\[
J=(J_o^{1/2}+\Delta_{\hat J})^2
~~,~~
\Delta_{\hat J}\in\left\{~\Delta=\Delta^T~:~\bmat{cc} \1 &\Delta \emat\bmat{cc} X^*&Y^*\\Y^*&Z^*\emat\bmat{c} \1\\\Delta\emat\preceq\0~\right\}.
\]
LFT modeling with respect to this newly defined uncertainty is rather simple following the same lines as the first method and gives
\[
\bmat{cc}J & J^{1/2}L\\L^TJ^{1/2}&\1\emat
=
\bmat{cc} \Delta_{\hat J}&\0\\\0&\Delta_{\hat J}\emat
\star
\bmat{cc|cc}
\0&\1&J_o^{1/2}&L\\
\0&\0&\1&\0\\\hline
\1&J_o^{1/2}&J_o&J_o^{1/2}L\\
\0&L^T&L^TJ_o^{1/2}&\1
\emat.
\]
The LFT built in this way has two remarkable features: \emph{i)} to the best of our knowledge, it is the first time that the modeling involves an uncertain matrix that is constrained to be symmetric, \emph{ii)} this matrix is repeated twice $\bmat{cc} \Delta_{\hat J}&\0\\\0&\Delta_{\hat J}\emat=\Delta_{\hat J}\otimes\1_2$. To build LMI type results for such uncertainties one needs to build some $DG$-scaling like result \cite{Fan91}. That is, to characterize via linear matrix equalities and inequalities the matrices $\Theta_{\hat J}$ that satisfy
\[\bary{l}
\bmat{cc} \1 & \Delta_{\hat J}\otimes\1_2 \emat \Theta_{\hat J} \bmat{c} \1 \\ \Delta_{\hat J}\otimes\1_2 \emat\preceq\0
\\
\forall \Delta_{\hat J}\in
\left\{~\Delta=\Delta^T~:~\bmat{cc} \1 &\Delta \emat\bmat{cc} X^*&Y^*\\Y^*&Z^*\emat\bmat{c} \1\\\Delta\emat\preceq\0~\right\}.
\eary\]
A choice of such matrices $\Theta_{\hat J}$ is natural generalization of the well-known $DG$-scalings that work for scalar repeated uncertainties
\[
\Theta_{\hat J}=
\bmat{cc}
X^*\otimes D & Y^*\otimes D+\1\otimes G
\\
Y^*\otimes D-\1\otimes G & Z^*\otimes D
\emat
~~:~~
\bary{l}
D=D^T\succ\0\in\mathbb{R}^{2\times2}
\\
G=-G^T\in\mathbb{R}^{2\times2}.
\eary
\]
The proof of this fact is trivial: in the formula following the $G$ dependent terms cancel one another thanks to the fact that $\Delta$ is symmetric and remains only
\[
\bmat{cc} \1 & \Delta\otimes\1_2 \emat \Theta \bmat{c} \1 \\ \Delta\otimes\1_2 \emat
=
D\otimes \left(\bmat{cc} \1 &\Delta \emat\bmat{cc} X^*&Y^*\\Y^*&Z^*\emat\bmat{c} \1\\\Delta\emat\right),
\]
which is negative semi-definite because it is the result of a Kronecker product of positive definite matrix and a negative semi-definite matrix.

\subsection{LFT modeling of the uncertain system}

Based on the described modeling of uncertain matrices discussed in the previous section and with some rather trivial additional manipulations---independent from  the choice for modeling the inertia $J$---the system dynamics can be converted to the following descriptor state-space form
\begin{equation}\label{eq: first state space model}
\left(\Delta_E\star\bmat{cc} E_d&E_c\\E_b&E_a\emat\right)\dot X
=
\left(\Delta_A\star\bmat{cc} A_d&A_c\\A_b&A_a\emat\right) X
+
Bu,
\end{equation}
where $X=\bvec{cccc}\dot\theta^T&\dot\eta^T&\theta^T&\eta^T\evec^T$ is the state of the satellite including its flexible modes;
$\Delta_A={\sf diag}\bmat{ccc} \delta_\Omega & \delta_Z & \delta_\Omega\emat$;
$\Delta_E={\sf diag}\bmat{cccc}  \delta_{J_1} & \delta_{J_1} & \delta_{\hat J_2} & \delta_{\hat J_2} \emat$ or
$\Delta_E=\Delta_{\hat J}\otimes \1_2$ depending on the choice of modeling of inertia;
$E$ and $A$ matrices are build accordingly. Taking the inverse of the left-hand side of \eqref{eq: first state space model} this formula allows to build a usual state-space model
\small
\[
\dot X=
\left(
{\sf diag}\bmat{c} \Delta_E \\\Delta_A\emat
\star
\bmat{cc|cc}
E_d-E_cE_a^{-1}E_b&-E_cE_a^{-1}A_b&-E_cE_a^{-1}A_a&-E_cE_a^{-1}B\\
\0 & A_d & A_c&\0\\\hline
E_a^{-1}E_b&E_a^{-1}A_b&E_a^{-1}A_a&E_a^{-1}B
\emat
\right)
\bvec{c} X\\u\evec,
\]
\normalsize
which is the same as the following linear system
\small
\[
\left\{\bary{r@{\,=\,}r@{\,+\,}r@{\,+\,}r}
\dot X
& E_a^{-1}A_aX
&\bmat{cc} E_a^{-1}E_b&E_a^{-1}A_b \emat w_\Delta
&E_a^{-1}Bu
\\
z_\Delta
&\bmat{c} -E_cE_a^{-1}A_a\\A_c\emat X
&\bmat{cc} E_d-E_cE_a^{-1}E_b&-E_cE_a^{-1}A_b\\\0 & A_d\emat w_\Delta
&\bmat{c} -E_cE_a^{-1}B \\\0\emat u
\eary\right.,
\]
\normalsize
in feedback loop with the uncertainty $w_\Delta={\sf diag}\bmat{c} \Delta_E \\\Delta_A\emat z_\Delta$.
Such system with feedback uncertainties can be easily defined in the \RR~toolbox. A dedicated function has been developed that outputs this model. The output is of the following type
\[
\left\{\bary{r@{\,=\,}r@{\,+\,}r@{\,+\,}r}
\dot X&AX&B_\Delta w_\Delta&B_u u
\\
z_\Delta&C_\Delta X&D_{\Delta\Delta}w_\Delta&D_{\Delta u}u
\\
y&C_yX&D_{y\Delta}w_\Delta&D_{y u}u
\eary\right.
~~,~~
w_\Delta=\Delta z_\Delta.
\]

\subsection{Reduced size variations of the uncertain model}\label{sec: reduced size model}

\begin{figure}
\myfig{0.5 \columnwidth}{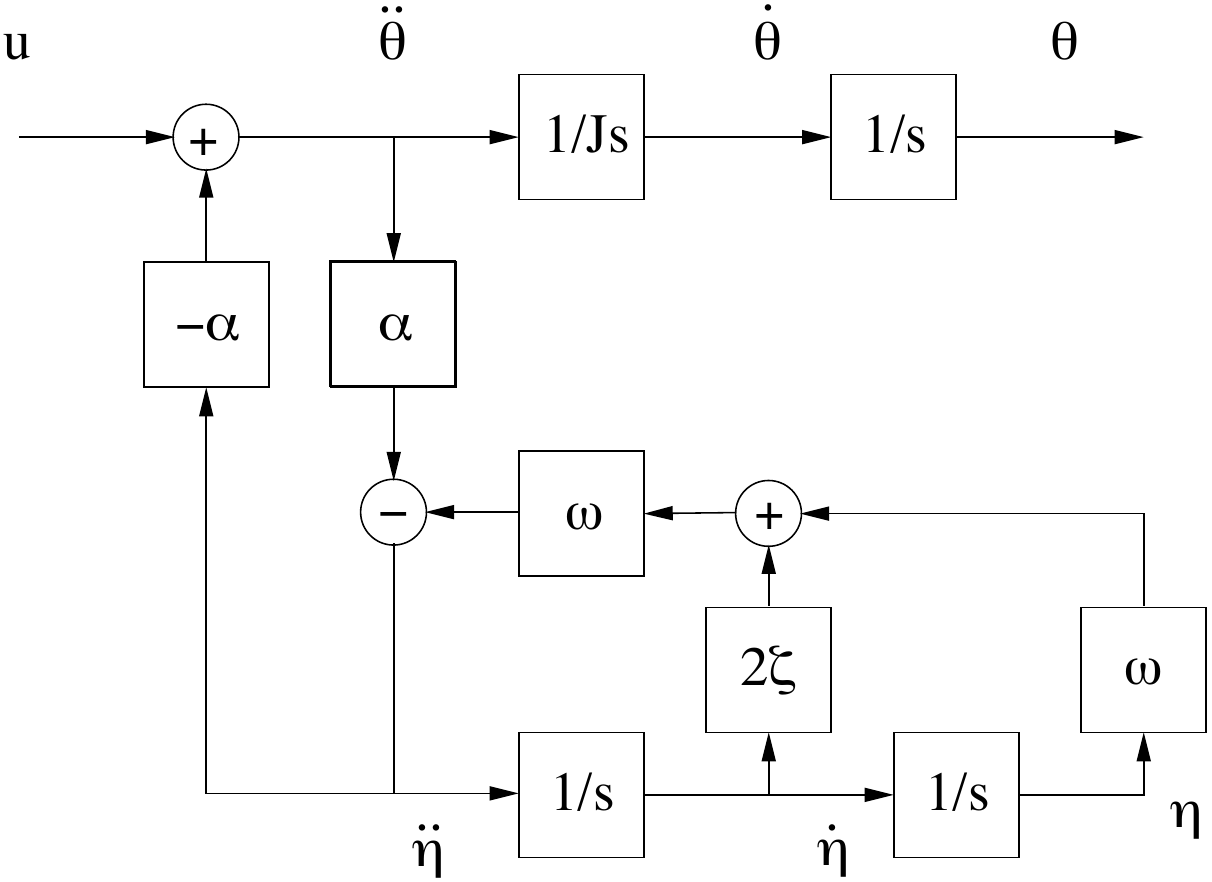}
\caption{Block diagram of one axis model with one flexible mode.}
\label{1DmodelOL}
\end{figure}

In order to test methods with respect to dimensions of the problem to solve (both in terms of order of the systems and in terms of size of the uncertainty block) several variants have been coded. The variations are threefold:
\begin{itemize}
\item[(a)] Select only one or two of the three axes. This of course reduces the number of states describing the satellite attitude. Moreover, in the case when only one axis is considered, the torsion and bending effects of the flexible modes can be combined. It produces models with twice less flexible modes states and twice smaller matrices $\Delta_A$.
\item[(b)] Select only some of the appendices. One can (virtually of course) remove any of the appendices. It produces models with reduced number of flexible modes and smaller matrices $\Delta_A$.
\item[(c)] Impose that all appendices have the same frequency and damping characteristics, $\omega_i=\omega$ and $\zeta_i=\zeta$. In such case, the number of flexible modes can be reduced to only three modes (one per axis) that are the projections of all bending and torsion modes on the attitude axes.
\end{itemize}
The simplest and rather realistic models amount to assuming (a) zero cross influence between satellite axes  and (c) that all appendices have exactly identical characteristics. Such assumptions reduce the study to three fourth-order models, one per angular axis. Each of these models ($i=1,2,3$) are described by two scalar equations
\begin{equation}\label{eq: simplest model}
\left\{\bary{l}
J_{ii}\ddot\theta_i+\sqrt{J_{ii}}l_i\ddot \eta_i=u_i\\
\sqrt{J_{ii}}l_i\ddot\theta_i+\ddot\eta_i +2\zeta\omega\dot\eta_i+\omega^2\eta_i=0
\eary\right.,
\end{equation}
and illustrated on Figure \ref{1DmodelOL} (where $\alpha=\sqrt{J_{ii}}l_i$). Corresponding LFT models have a $5\times5$ uncertain matrix where scalar uncertainties on $J_{ii}$ appear twice, scalar uncertainties on $\omega$ appear twice and scalar uncertainties on $\zeta$ appear once.

\section{State-feedback design model}
The control design problem is to build a control that ensures the following performances
\begin{itemize}
\item[(i)] As small as possible pointing error. To this end, the control should contain an integrator to improve the low frequency disturbing torques rejection.
\item[(ii)] Avoid saturation of the reaction wheel actuators. These actuators have the following nonlinear model
\[
u=sH(s) ~{\sf sat}_W\left( \frac{1}{s} {\sf sat}_T\left( u_c \right)\right),
\]
where $u_c$ stands for the torque control input computed by the controller and $u$ is the actual torque applied by the reaction wheel. ${\sf sat}_T$ is a saturation on the torque to be applied which is of $5\times 10^{-3}$Nm. It is in general not critical and can be neglected. The term $\frac{1}{s}$ is an integrator that outputs the reaction wheel angular momentum. This angular momentum is saturated (${\sf sat}_W$) with saturation level of $0.12$ Nms. This saturation is critical: when it occurs the system is no more actuated and is open-loop unstable. Finally, $sH(s)$ is a transfer function describing the dynamics of the reaction wheel.
\item[(iii)] Other specifications such as noise rejection, robustness to time-delays in the control, etc. as discussed in \cite{pit:arz/rocond06}. 
\end{itemize}

\begin{figure}
\myfig{0.5 \columnwidth}{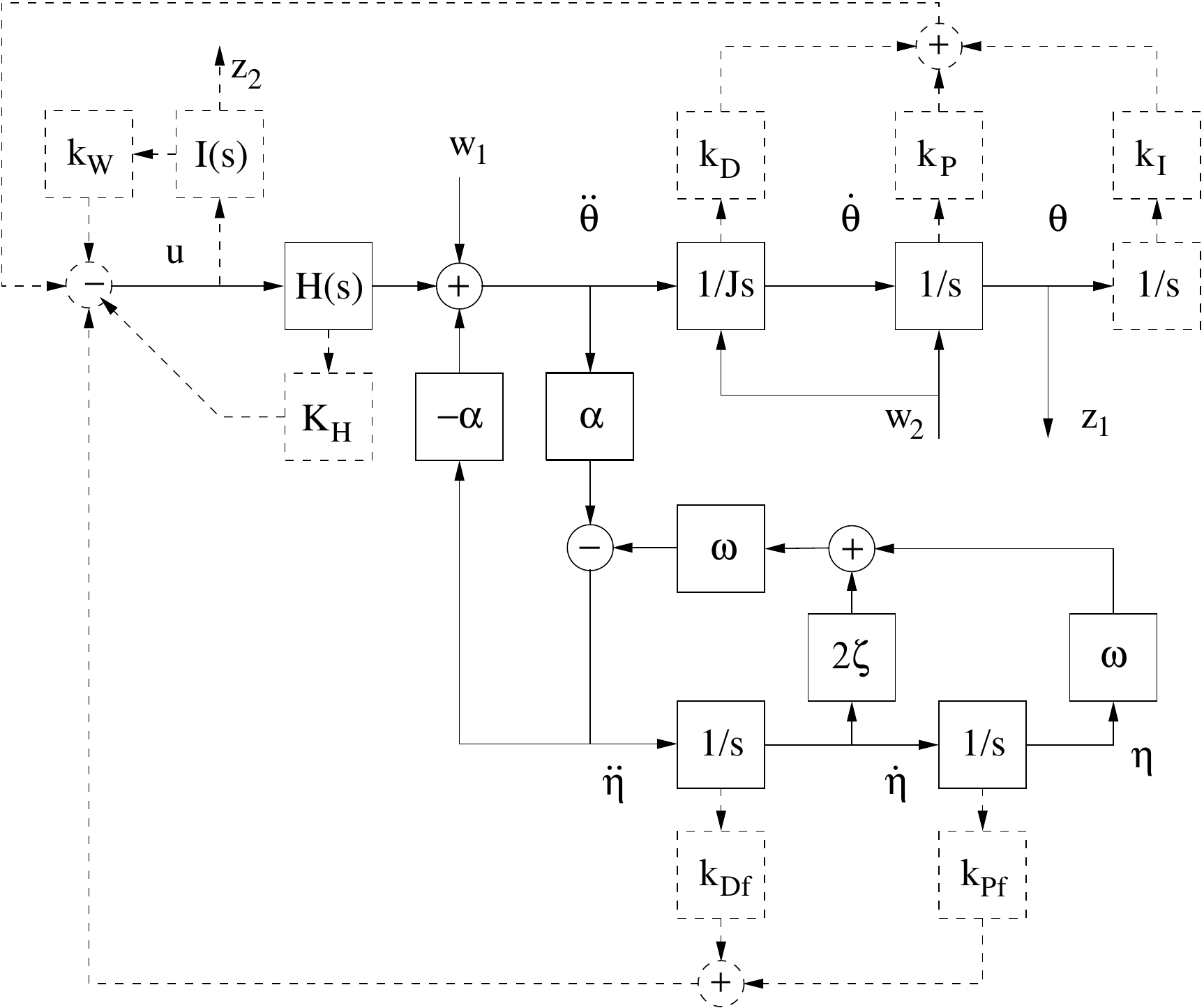}
\caption{Block diagram of state-feedback design model.}
\label{1DmodelSFB}
\end{figure}

In order to take into account the two specifications (i) and (ii), we add to the model an integrator of the output and a pseudo integrator $I(s)=\frac{1}{s+0.001}$ of the input. We remark that an integrator in the input---instead of pseudo integrator---would result in instability since the states of integrator are not controllable in the formulation. These are represented with dotted lines on Figure \ref{1DmodelSFB}. The dotted lines indicate that these blocks are added by the designer and hence part of the control law. 

For that augmented model we search for a robust state-feedback control as illustrated in Figure \ref{1DmodelSFB}. The dotted lines represent the state-feedback with eight gains. $k_P$, $k_I$, $k_D$ are the feedback gains with respect to the angular error $\theta$, the integral of it and its derivative respectively. $k_{Pf}$ and $k_{Df}$ are the gains on the angular position of the flexible mode $\eta$ and on its derivative respectively. $k_W$ is the gain on the state of the pseudo-integrator that models the reaction wheel speed. $K_H\in\mathbb R^{1\times 2}$ is the gain on the states of the reaction wheels. The aim of the control is to minimize the peak of $z_2$ (the reaction wheel speed) especially when the satellite starts from large non-zero angle and angular rate initial conditions which are represented as input signals $w_2$. We assume a maximal $\pm0.08$deg/s angular rate initial deviation and $\pm15$deg angular initial deviation. Simultaneously, the control should minimize the effect of unknown input perturbations on the system precision, that is to minimize the transfer for $w_1$ to $z_1$.

The design of such state-feedback controller is possible using \RR~toolbox \cite{rromulocconf,romulocconf}. In particular a function named \texttt{demeterPerformance} is developed to generate models required for controller design. The following lines of codes define three models being
\begin{enumerate}
\item
The augmented model with integerator on the output, reaction wheel model and pseudo-integrator of the input.
\item
Model with $w_1/z_1$ performance input output. 
\item
Model with $w_2/z_2$ performance input output. 
\end{enumerate}

\begin{verbatim}
usysIW=demeterPerformance(ConsideredAxis,ConsideredAppendices,...
    model_type,uncertainty_type, rwheels,0);
usysIW1=demeterPerformance(ConsideredAxis,ConsideredAppendices,...
    model_type,uncertainty_type, rwheels,1);
usysIW2=demeterPerformance(ConsideredAxis,ConsideredAppendices,...
    model_type,uncertainty_type, rwheels,2);
\end{verbatim}
Next, we briefly explain various arguments of the \texttt{demeterPerformance} function.\\ The parameters \texttt{ConsideredAxis} and \texttt{ConsideredAppendices} define the number of axes and appendices used in the model respectively. If \texttt{model\_type=2},  all flexible modes  have the same frequency and damping characteristic with the same uncertain parameters but, if  \texttt{model\_type=1}, uncertain parameters are allowed  to be independent for different appendices. If \texttt{uncertainty\_type=1}, all uncertainties are norm-bounded scalars, if \texttt{uncertainty\_type=2}, all uncertainties are scalars in intervals and if \texttt{uncertainty\_type=3}, uncertainties on inertia are norm-bounded deterministic, others are uniformly distributed in intervals. If \texttt{rwheels=1}, the reaction wheels are included in the model and if \texttt{rwheels=0}, the model does not include reaction wheels dynamic. 

Let $N_a$ be the number of considered axes and $N_f$ be the number of appendices. The satellite dynamics involve $2*N_a+4*N_f$ states to which one adds actuator models and $N_a$ integrators of the control law. In case \texttt{model\_type=1} (all appendices have different characteristics) the satellite dynamics involve $N_a(N_a+1)/2+2*N_f$ scalar uncertainties. If \texttt{model\_type=2} (all appendices have identical characteristics) the satellite dynamics involve $N_a(N_a+1)/2+2$ scalar uncertainties. A special case is when $N_a=1$ and all appendices are considered identical. In such case the satellite dynamics involve only $4$ states and $3$ uncertainties, see \eqref{eq: simplest model}.

In \RR~there are two approaches to design the robust state feedback controller. The first approach is based on deterministic multiobjective methods in which the performance specifications are enforced to hold for the entire set of uncertainty. The second paradigm is probabilistic and randomized methods in which the design specifications (including stability) are enforced to hold up to a probability level. In the next two subsections, we study the two mentioned approaches in state feedback design.
\section{Controller design}
\subsection{Deterministic approach}\label{sec: deterministic design}
In \RR~ the deterministic state-feedback design LMI problem is defined as
\begin{verbatim}
quiz=ctrpb('state-feedback','unique')...
  +1*hinfty(usysIW1)...
  +100*i2p(usysIW2)...
  +dstability(usysIW,region('plane',-1e-4))...
  +dstability(usysIW,region('plane',-10,pi));
\end{verbatim}
The LMI problem built in this way is based on quadratic stability type results with Lyapunov shaping paradigm \cite{scherer_multiobjective_1997}, that is, a unique Lyapunov matrix is used for assessing all four specified performances and for all values of uncertainties. The four specifications are: the $\mathcal{H}_\infty$ performance with respect to the input/outputs $w_1/z_1$; the impulse-to-peak performance with respect to the input/outputs $w_2/z_2$ (which is equivalent to looking at peak response to the initial conditions); the pole location performance such that all closed-loop poles should have real part smaller than $-1\times10^{-4}$ and greater than $-10$ (influences the rapidity of the time response).
The LMI problem is solved in \RR~using the following commands that returns the state-feedback gain
\begin{verbatim}
Ksf_det=solvesdp(quiz,sdpsettings('verbose',1,'solver','mosek'));
\end{verbatim}

\subsection{Probabilistic Design }\label{sec: probabilistic design}
There are two paradigms in probabilistic techniques for controller design. The first approach is non-sequential, in which a sampled version of the original problem is solved in one shot. The scenario approach \cite{calafiore_uncertain_2004,calafiore_scenario_2006} is a non-sequential approach for solving uncertain convex problems. The main idea in this approach is to reformulate a semi-infinite convex optimization problem as a sampled convex optimization problem subject to a finite number of random constraints extracted form the uncertainty set.
The second class of probabilistic  design algorithms are sequential methods, in which at each iteration, a candidate solution is constructed---based on gradient \cite{polyak_probabilistic_2001}, ellipsoid \cite{kanev_ellipsoid_2003}, cutting plane \cite{dabbene_randomized_2010} or sampling based technique \cite{chamanbaz_tac_sequential_2013}---and its robustness is verified through a sequential probabilistic validation algorithm \cite{alamo2015randomized}. In \RR, the scenario approach and sequential algorithms based on gradient update rule \cite{polyak_probabilistic_2001} and the sequential approach presented in \cite{chamanbaz_tac_sequential_2013} are used to solve uncertain state-feedback design problem. 
A controller addressing the same performance requirements as in the deterministic case can be formulated and solved using the sequential algorithm  \cite{chamanbaz_sequential_2013,chamanbaz_tac_sequential_2013} 
\begin{verbatim}
quiz = ctrpb('state-feedback','rand')...
  +1*hinfty(usysIW1)...
  +100*i2p(usysIW2)...
  +dstability(usysIW,region('plane',-1e-4))...
  +dstability(usysIW,region('plane',-10,pi));
opts=randsettings('epsilon',0.1,'delta',1e-9,...
'method','sequential','sdpopts',...
sdpsettings('verbose',0,'solver','mosek'));
Ksf_prob=solvesdp(quiz,opts);
\end{verbatim}
The parameters \texttt{epsilon} and \texttt{delta} defined in the \texttt{randsettings} function are the required accuracy and confidence levels of the solution.
In words, the probability that the solution does not satisfy constraints is smaller than \texttt{epsilon} and this statement holds with probability at least \texttt{1-delta}. We refer to \cite{calafiore_research_2011,tempo_randomized_2012} for the exact definition of accuracy and confidence levels.
We remark that one can solve the same problem using the scenario approach \cite{calafiore_uncertain_2004,calafiore_scenario_2006} by changing \texttt{`sequential'} to \texttt{`scenario'} in the code. 

\section{Closed-loop analysis of the state-feedback law}\label{sec: analysis}
An important feature of \RR~is to provide in a unified framework different available tools for analyzing the robust performance of a uncertain closed-loop systems. In particular, a user can check if several  performance criteria, as for instance $\mathcal{H}_2$ and $\mathcal{H}_\infty$ norms, impulse-to-peak response, pole location, etc., hold either robustly or with a guaranteed level of probability. 
Similar to design techniques, analysis can be performed either in a deterministic setting or through randomized algorithms resulting in a probabilistic estimate of robust performance. 
\subsection{Deterministic analysis}
The deterministic analysis methods implemented in \RR~ are based on Lyapunov-type certificates. In particular, it can be based on either parameter dependent Lyapunov function \cite{Ebi14-SVbook,Iwa01,Pea07a} or a common Lyapunov function \cite{scherer_multiobjective_1997}.
An upper bound of the closed-loop $\mathcal{H}_\infty$ norm for the transfer $z_1/w_1$ can be computed using parameter-dependent Lyapunov matrices as follows 
\begin{verbatim}
usysIW1cl=sfeedback(usysIW1,Ksf_det); 
quiz = ctrpb('analysis', 'PDLF')+hinfty(usysIW1cl); 
solvesdp(quiz,sdpopts);
\end{verbatim}

\subsection{Probabilistic analysis} \label{sec: randomized analysis}
The probabilistic analysis is based on a Monte Carlo algorithm in which a number of random samples are extracted from the set of uncertainty and the performance index is measured only at the extracted samples. There are two probabilistic analysis algorithms: \textit{1)} Worst-case performance estimation in which an estimate of the worst-case performance is defined as the worst-case performance among all extracted samples. The sample size in this case is defined by a log-over-log bound \cite{TeBaDa:97}. \textit{2)} Randomized performance verification where the objective is to estimate  the probability of a given level of performance being satisfied, for instance estimating the probability of instability or the probability that the $\mathcal{H}_\infty$ norm of the system is below a given level. The number of samples in this case is defined by the Chernoff bound \cite{Chernoff_1952}.  The next command computes the wost-case $\mathcal{H}_\infty$ norm of the closed-loop system \texttt{usysIW1cl} using a randomized worst-case performance estimation algorithm.
\begin{verbatim}
quiz = ctrpb('analysis', 'rand')+hinfty(usysIW1cl); 
opts=randsettings('epsilon',1e-1,'delta',1e-6);
solvesdp(quiz,opts);
\end{verbatim}

\begin{table*}[!t]
\begin{center}
\scalebox{0.8}{
\begin{tabular}{cccccccccccccc}
\rotatebox{90}{\texttt{ConsideredAxis}} & \rotatebox{90}{\texttt{ConsideredAppendices}} & \rotatebox{90}{\texttt{model\_type}} & \rotatebox{90}{\texttt{uncertainty\_type}} & \rotatebox{90}{\texttt{rwheels}} &   & \multicolumn{2}{c}{\emph{Design}} & \multicolumn{4}{c}{\emph{Analysis}} & \tabularnewline
\midrule

 & & & & & & & & \multicolumn{2}{c}{\emph{Det}} & \multicolumn{2}{c}{\emph{Prob}} &  \emph{Complexity(s)} \tabularnewline

\midrule
& & & &  & \emph{Design} & \emph{impulse} & \emph{Infinity} & \emph{impulse} & \emph{Infinity} & \emph{impulse} & \emph{Infinity } & \tabularnewline
& & & &  & \emph{Method} & \emph{to Peak} & \emph{Norm} & \emph{to Peak} & \emph{Norm} & \emph{to Peak} & \emph{Norm} & \tabularnewline
\midrule
\midrule
1 & 1 & 1 & 1 & 1 & Prob & 22.3 & 2.9 & 0.36 & 1.5 & 0.13 & 1.01 & 160  \tabularnewline
1 & 1 & 1 & 1 & 1 &  Det & 22.3 & 4.7 & 0.41 & 1.5 & 0.16 & 1.16 & 1 \tabularnewline
\midrule
1 & 1,2 & 1 & 1 & 1 &  Prob & Inf & Inf & NA & NA & NA & NA & NA \tabularnewline
1 & 1,2 & 1 & 1 & 1 &  Det & Inf & Inf & NA & NA & NA & NA & NA \tabularnewline
\midrule
1 & 1,2,3,4 & 2 & 1 & 1 &  Prob & 22.5 & 3 & 0.42 & 1.3 & 0.14 & 0.84 & 520 \tabularnewline
1 & 1,2,3,4 & 2 & 1 & 1 &  Det & 22.5 & 3 & 0.43 & 1.3 & 0.13 & 0.99 & 1.3 \tabularnewline
\midrule

1,2 & 1,2,3,4 & 2 & 1 & 1 &  Prob & Inf & Inf & NA & NA & NA & NA & NA \tabularnewline
1,2 & 1,2,3,4 & 2 & 1 & 1 &  Det & Inf & Inf & NA & NA & NA & NA & NA \tabularnewline
\midrule
1,2 & 1,2 & 2 & 1 & 1 &  Prob & 22.4 & 2.8 & 0.67 & Inf & 0.2 & 0.06 & 2215 \tabularnewline
1,2 & 1,2 & 2 & 1 & 1 &  Det & 22.7 & 5 & Inf & Inf & 0.16 & 0.5 & 142 \tabularnewline
\midrule
1,2 & 1,2 & 2 & 2 & 1 &  Prob & 22.46 & 2.69 & 0.7 & 1.38 & 0.19 & 0.08 & 1750 \tabularnewline
1,2 & 1,2 & 2 & 2 & 1 &  Det & 22.6 & 4.24 & 0.75 & 1.03& 0.19 & 0.14 & 46 \tabularnewline
\midrule
1,2,3 & 1,2 & 2 & 2 & 1 &  Prob & 22.5 & 3.3 &  Inf & Inf & 0.23 & 0.66 & 16387 \tabularnewline
1,2,3 & 1,2 & 2 & 2 & 1 &  Det & 22.7 & 8.1 & Inf & Inf & 0.2 & 1.34 & 14111 \tabularnewline

\bottomrule
\end{tabular}
}
\end{center}
\caption{Simulation results for various probabilistic and deterministic controllers designed using \RR~for the DEMETER model. ``Inf'' indicates the cases where the optimization problem is infeasible; ``NA'' also refers to Not Applicable.}
\label{tab: simulation results}
\end{table*}

\section{Numerical tests}\label{sec: comparison}
In this section, we compare probabilistic and deterministic approaches in terms of performance and complexity. To this end, we generate a number of DEMETER models---based on the discussion of subsection \ref{sec: reduced size model} by changing parameters \texttt{ConsideredAxis}, \texttt{ConsideredAppendices}, \texttt{model\_type}, \texttt{uncertainty\_type} and \texttt{rwheels}---and design different deterministic and probabilistic controllers. Next, the performance of designed controllers is measured using deterministic and probabilistic analysis methods of section \ref{sec: analysis} to quantify the level of conservatism associated with different design approaches. The result of these numerical tests is reported in Table \ref{tab: simulation results} where we consider different number of axes, appendices and different model and uncertainty types and design probabilistic and deterministic controllers for the generated models. 
The probabilistic controller is designed using the scenario approach and probabilistic accuracy \texttt{epsilon} and confidence \texttt{delta} levels are set to $0.1$ and $10^{-9}$ respectively.   
In most cases---as expected---the probabilistic controller achieves less conservative performance levels in handling various uncertainties. 
In terms of computational complexity, the deterministic approach is less computationally demanding for the case that all  uncertainties are considered to be norm bounded.
However, if we require  uncertainties to be defined in intervals (and hence in polytopes), the computational complexity associated with the deterministic approach increases significantly. For such uncertainties, \RR~applies a vertex-separator result as proposed in \cite{Iwa98}. At the difference of highly sparse DG-scaling type separators with few constraints that are build in case of norm-bounded uncertainties, the vertex-separator is known to be less conservative but with increased number of decision variables (full matrices) and increased number of constraints (one for each vertex, and the number of vertices is $2^N$ where $N$ is the number of uncertain parameters).  
We remark that in some problem instances of Table \ref{tab: simulation results} that the optimization problem---for controller design---is infeasible there does not exists a ``robust'' state-feedback controller satisfying all required specifications and the optimization problem becomes infeasible even for large probabilistic accuracy \texttt{epsilon} and confidence \texttt{delta} levels. 
\begin{figure}[t]
\centering
\includegraphics[width=0.6\columnwidth]{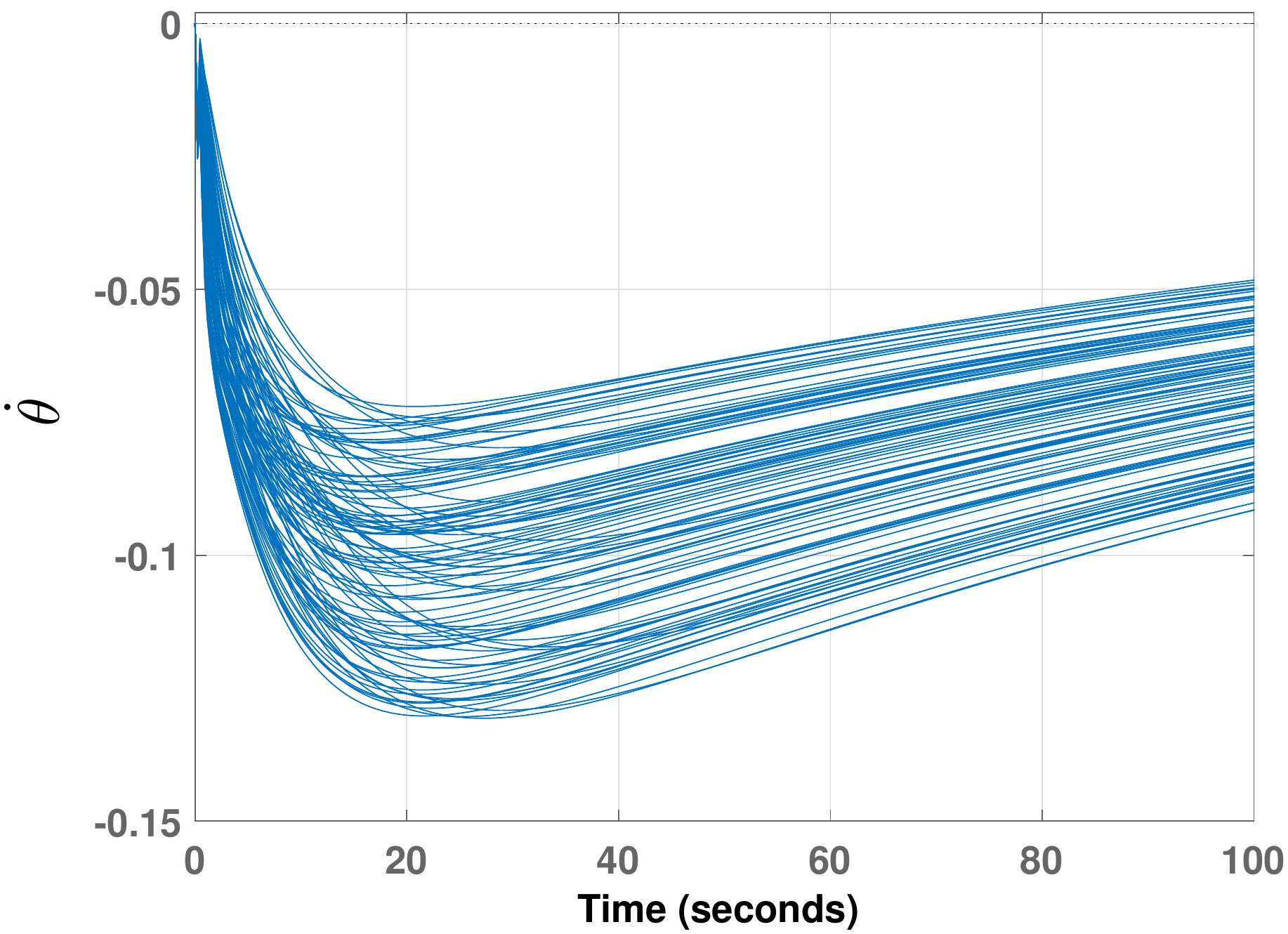}
\caption{ Impulse response of 100 randomly generated closed-loop systems from $w_2$ to $z_2$ with the controller designed in second row of Table \ref{tab: simulation results}.  }
\label{fig: Impulse response}
\end{figure}
\begin{figure}[t]
\centering
\includegraphics[width=0.6\columnwidth]{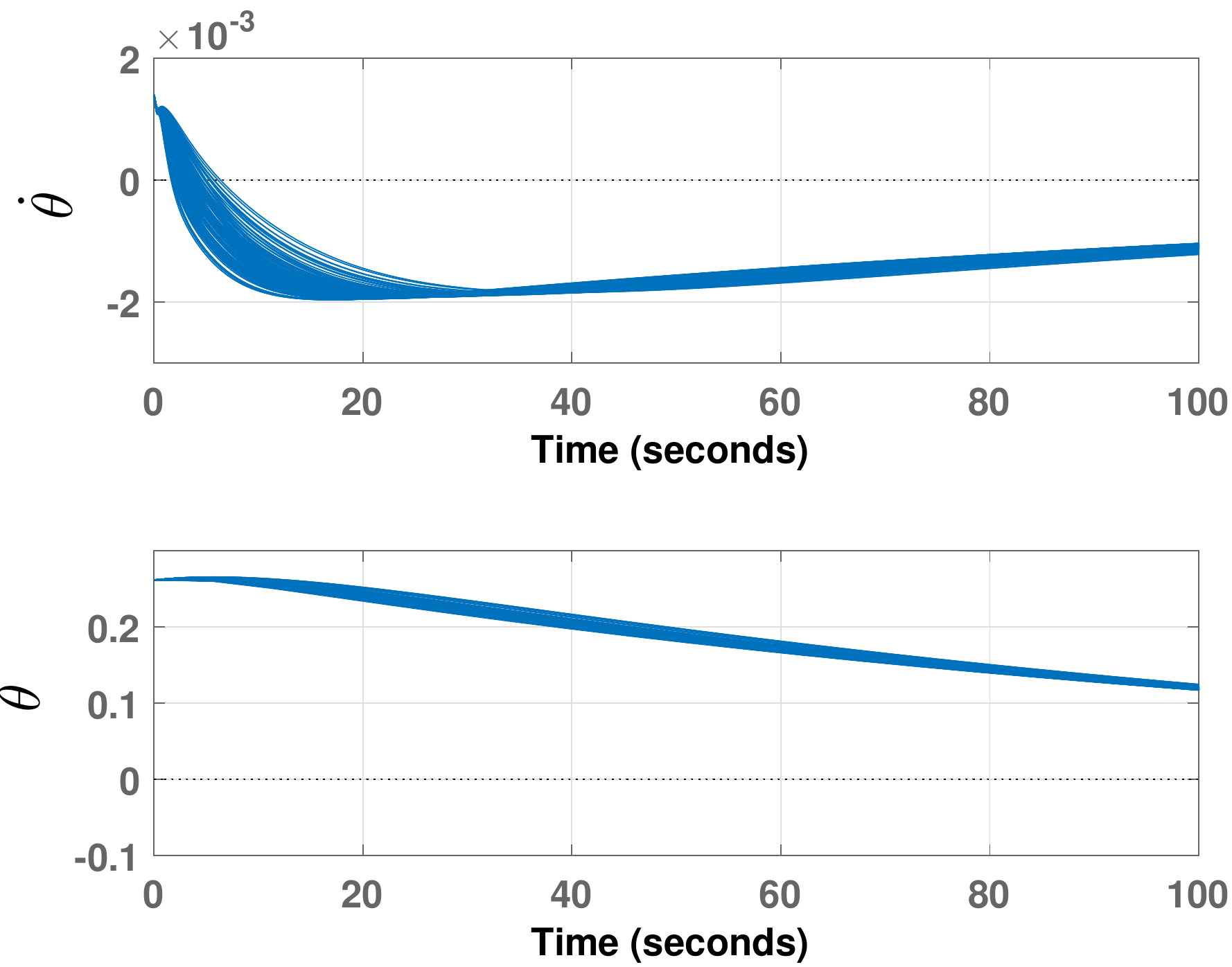}
\caption{ Time trajectories of satellite angular rate $\dot{\theta}$ (top figure) and angular deviation $\theta$ (bottom figure) for 100 randomly generated closed-loop systems from $w_2$ to $z_2$ with the controller designed in second row of Table \ref{tab: simulation results}. }
\label{fig: Time trajectory}
\end{figure}

To further validate our design, a posteriori analysis using Monte-Carlo simulation is carried on for the controller designed in the second row of Table \ref{tab: simulation results}. To do so, we extracted 100 random samples from the uncertainty set, closed the loop for each of them and measured the impulse response---from $w_2$ to $z_2$---of each sampled closed-loop system. Figure \ref{fig: Impulse response} shows the result of this simulation. Figure \ref{fig: Time trajectory} also 
demonstrates the time trajectories of the angular rate $\dot{\theta}$ and angular deviation $\theta$ of the satellite for the same sampled closed-loop systems.  One can see that $\dot\theta$ starts from the initial condition $0.08\pi/180=1.4\times 10^{-3}$rad/s and $\theta$ starts from $15\pi/180=0.262$rad. This is considered as the worst case initial configuration. It is such that the pointing error $\theta$ tends to increase at the start due to the positive angular rate.

An interesting feature of randomized methods is that the computational complexity does not depend on the number of uncertain parameters. This feature is known as breaking the curse of dimensionality. Therefore, increasing the number of uncertain parameters does not influence the complexity of solving state-feedback problem using randomized methods. On the other hand, stability and performance  achieved using the controller designed by this approach is not guaranteed to hold for the entire set of uncertainties. That is, there might exist a subset of the uncertain set---although with very small probability measure---for which the guaranteed performance level is violated.

It is noted that the designed controllers in this paper are of state-feedback type, requiring all the states to be available for feedback. This requirement is not realistic in practice. In fact, in practice,  sensors report $\theta, \dot{\theta}$ and  $\int\theta$. Observers are needed for flexible modes $\eta, \dot{\eta}$. Therefore, an observer can be designed using the approach presented in  \cite{peaucelle2014lmi} in order to estimate the states of the system and then use the state-feedback controller formulated in this paper to control the DEMETER satellite.

\section{Conclusions}

This paper shows how the features of the recently release Matlab toolbox \RR~can be exploited to perform both deterministic and probabilistic analysis and design of systems in the presence of uncertainty. The potentialities of \RR~ are illustrated on the DEMETER satellite benchmark. The performed numerical simulations are fully reproducible, since both the DEMETER model and the \RR~toolbox are freely downloadable at \url{http://projects.laas.fr/OLOCEP/rromuloc/}.

\section*{Acknowledgments}
Acknowledgments to all those who contributed to {\RR} in many different ways:
D. Arzelier,  A. Bortott,  G. Calafiore,  G. Chevarria,  E. Gryazina,  B. Polyak,   P. Shcherbakov,   M. Sevin,  P. Spiesser,  and A. Tremba.

\bibliographystyle{plain}
\bibliography{ref}

\begin{thebibliography}{10}

\bibitem{alamo2015randomized}
T.~Alamo, R.~Tempo, A.~Luque, and D.R. Ramirez.
\newblock Randomized methods for design of uncertain systems: Sample complexity
  and sequential algorithms.
\newblock {\em Automatica}, 52:160--172, 2015.

\bibitem{calafiore_uncertain_2004}
{G.C.} Calafiore and {M.C.} Campi.
\newblock Uncertain convex programs: randomized solutions and confidence
  levels.
\newblock {\em Mathematical Programming}, 102:25--46, 2004.

\bibitem{calafiore_scenario_2006}
{G.C.} Calafiore and {M.C.} Campi.
\newblock The scenario approach to robust control design.
\newblock {\em {IEEE} Transactions on Automatic Control}, 51:742--753, 2006.

\bibitem{calafiore_research_2011}
{G.C.} Calafiore, F.~Dabbene, and R.~Tempo.
\newblock Research on probabilistic methods for control system design.
\newblock {\em Automatica}, 47:1279--1293, 2011.

\bibitem{rromulocconf}
M.~Chamanbaz, F.~Dabbene, D.~Peaucelle, and R.~Tempo.
\newblock {R-RoMulOC}: a unified tool for randomized and robust multiobjective
  control.
\newblock In {\em 8th IFAC Symposium on Robust Control Design}, Bratislava,
  July 2015.

\bibitem{chamanbaz_sequential_2013}
M.~Chamanbaz, F.~Dabbene, R.~Tempo, V.~Venkataramanan, and Q-G. Wang.
\newblock Sequential randomized algorithms for sampled convex optimization.
\newblock In {\em Proc. {IEEE} Multi-Conference on Systems and Control}, pages
  182--187, Hyderabad, India, 2013.

\bibitem{chamanbaz_tac_sequential_2013}
M.~Chamanbaz, F.~Dabbene, R.~Tempo, V.~Venkataramanan, and Q-G. Wang.
\newblock Sequential randomized algorithms for convex optimization in the
  presence of uncertainty.
\newblock {\em IEEE Transactions on Automatic Control}, 61:2565--2571, 2016.

\bibitem{Chernoff_1952}
H.~Chernoff.
\newblock A measure of asymptotic efficiency for tests of a hypothesis based on
  the sum of observations.
\newblock {\em The Annals of Mathematical Statistics}, 23:493--507, 1952.

\bibitem{dabbene_randomized_2010}
F.~Dabbene, P.~S. Shcherbakov, and B.~T. Polyak.
\newblock A randomized cutting plane method with probabilistic geometric
  convergence.
\newblock {\em {SIAM} Journal on Optimization}, 20, 2010.

\bibitem{Ebi14-SVbook}
Y.~Ebihara, D.~Peaucelle, and D.~Arzelier.
\newblock {\em S-variable Approach to {LMI}-based Robust Control}.
\newblock Springer London, 2015.

\bibitem{Fan91}
M.~Fan, A.~Tits, and J.~Doyle.
\newblock Robustness in the presence of mixed parametric uncertainty and
  unmodelled dynamics.
\newblock 36(1):25--38, January 1991.

\bibitem{Iwa98}
T.~Iwasaki and S.~Hara.
\newblock Well-posedness of feedback systems: Insights into exact robustness
  analysis and approximate computations.
\newblock {\em IEEE Trans. on Automat. Control}, 43(5):619--630, 1998.

\bibitem{Iwa01}
T.~Iwasaki and G.~Shibata.
\newblock {LPV} system analysis via quadratic separator for uncertain implicit
  systems.
\newblock {\em IEEE Transactions on Automatic Control}, 46:1195--1208, 2001.

\bibitem{kanev_ellipsoid_2003}
S.~Kanev, B.~De~Schutter, and M.~Verhaegen.
\newblock An ellipsoid algorithm for probabilistic robust controller design.
\newblock {\em Systems \& Control Letters}, 49:365--375, 2003.

\bibitem{romulocconf}
D.~Peaucelle and D.~Arzelier.
\newblock Robust multi-objective control toolbox.
\newblock In {\em Proc. of the {CACSD} Conference}, Munich, Germany, 2006.

\bibitem{Pea07a}
D.~Peaucelle, D.~Arzelier, D.~Henrion, and F.~Gouaisbaut.
\newblock Quadratic separation for feedback connection of an uncertain matrix
  and an implicit linear transformation.
\newblock {\em Automatica}, 43:795--804, 2007.

\bibitem{peaucelle2014lmi}
D.~Peaucelle and Y.~Ebihara.
\newblock {LMI} results for robust control design of observer-based
  controllers, the discrete-time case with polytopic uncertainties.
\newblock {\em In Proc. 19th IFAC world congress}, pages 6527--6532, 2014.

\bibitem{petersen_robust_2014}
I.R. Petersen and R.~Tempo.
\newblock Robust control of uncertain systems: Classical results and recent
  developments.
\newblock {\em Automatica}, 50:1315--1335, 2014.

\bibitem{pit:arz/rocond06}
C.~Pittet and D.~Arzelier.
\newblock Demeter: A benchmark for robust analysis and control of the attitude
  of flexible micro satellites.
\newblock In {\em Proc. 6th {IFAC} Symposium on Robust Control Design}, pages
  661--666, Toulouse, France, 2006.

\bibitem{polyak_probabilistic_2001}
B.~T. Polyak and R.~Tempo.
\newblock Probabilistic robust design with linear quadratic regulators.
\newblock {\em Systems \& Control Letters}, 43:343--353, 2001.

\bibitem{scherer_multiobjective_1997}
C.~Scherer, P.~Gahinet, and M.~Chilali.
\newblock Multiobjective output-feedback control via {LMI} optimization.
\newblock {\em {IEEE} Transactions on Automatic Control}, 42:896--911, 1997.

\bibitem{TeBaDa:97}
R.~Tempo, E.-W. Bai, and F.~Dabbene.
\newblock Probabilistic robustness analysis: {E}xplicit bounds for the minimum
  number of samples.
\newblock {\em Systems and Control Letters}, 30:237--242, 1997.

\bibitem{tempo_randomized_2012}
R.~Tempo, G.C. Calafiore, and F.~Dabbene.
\newblock {\em Randomized Algorithms for Analysis and Control of Uncertain
  Systems: With Applications}.
\newblock Springer, 2nd edition, 2013.

\bibitem{tremba_ract:_2008}
A.~Tremba, G.C. Calafiore, F.~Dabbene, E.~Gryazina, B.~Polyak, P.~Shcherbakov,
  and R.~Tempo.
\newblock {RACT:} randomized algorithms control toolbox for {MATLAB}.
\newblock In {\em Proc. 17th World Congress of {IFAC}, Seoul}, pages 390--395,
  2008.

\end{thebibliography}
\end{document}